\documentclass{article}
\usepackage{amsmath}
\usepackage{makecell}

\usepackage{tikz}
\usetikzlibrary{shapes.geometric, arrows, patterns}

\tikzstyle{subpos} = [rectangle, minimum width=3cm, minimum height=0.75cm, text centered, draw=black]
\tikzstyle{arrow} = [thick,-,>=stealth]
\usepackage{multirow,array}
\usepackage{cancel}
\usepackage[font=small,labelfont=bf]{caption}
\usepackage[font=small]{subcaption}

\usepackage{makeidx}         
\usepackage{graphicx}        
\usepackage[leftcaption]{sidecap}        
                    
\usepackage{multicol}        
\usepackage[bottom]{footmisc}
\usepackage[hidelinks]{hyperref}

\begin{document}

\title{Spatial Iterated Prisoner’s Dilemma as a
Transformation Semigroup}

\author{
  Isaiah Farahbakhsh
  \thanks{University of Waterloo, 200 University Ave W, Waterloo, ON N2L 3G1, Canada \hspace{6mm} email: infarahb@uwaterloo.ca}
  \and
  Chrystopher L. Nehaniv 
  \thanks{Algebraic Intelligence \& Computation Laboratory University of Waterloo, 200 University Ave W, Waterloo, ON N2L 3G1, Canada \hspace{3mm} email: chrystopher.nehaniv@uwaterloo.ca}
  \thanks{This work was supported by the Natural Sciences and Engineering Research Council of Canada (NSERC), funding ref. RGPIN-2019-04669, and the University of Waterloo}
}
\date{}
\maketitle
\begin{abstract}
    The prisoner’s dilemma (PD) is a game-theoretic model studied in a wide array of fields to understand the emergence of cooperation between rational self-interested agents. 
    In this work, we formulate a spatial iterated PD as a discrete-event dynamical system where agents play the
    game in each time-step and analyse it algebraically using Krohn-Rhodes algebraic automata theory 
     using a computational implementation of the holonomy decomposition of transformation semigroups. 
    In each iteration all players adopt the most profitable strategy in their immediate neighbourhood. Perturbations resetting the strategy of a given player provide additional generating events for the dynamics. \newline\indent Our initial study shows that the algebraic structure, including how natural subsystems 
    comprising permutation groups acting on the spatial distributions of strategies, arise in certain parameter regimes for the pay-off matrix, and are absent for other parameter regimes. Differences in the number of group levels in the holonomy decomposition (an upper bound for Krohn-Rhodes complexity) are revealed as more pools of reversibility appear when the temptation to defect is at an intermediate level. Algebraic structure uncovered by this analysis can be interpreted to shed light on the dynamics of the spatial iterated PD.
\end{abstract}

\section{Introduction}
\label{sec:1}

Krohn-Rhodes (KR) theory offers powerful tools for understanding discrete-event dynamical systems (e.g. \cite{nehaniv2015symmetry}). This theory decomposes any system whose dynamics can be represented as a transformation semigroup into a cascade of  permutation-group layers and identity-reset (flip-flop) layers using the wreath product \cite{krohn1965algebraic}. This yields a ``coarse-to-fine graining'' of both the system's state and its dynamical transformations.
The decomposition process can uncover subsystems represented by permutation groups which we call {\em pools of reversibility} or {\em natural subsystems} (see below). Algebraic structure uncovered by this analysis can be interpreted to shed light on the dynamics and complexity of many broad classes of discrete-event dynamical systems, including models found in the field of game theory.

The Prisoner's Dilemma (PD) is an extensively studied game which explores the problem of individual vs. collective profit in a simple 2-strategy model. The model is usually presented describing a situation where two partners-in-crime are imprisoned and unable to communicate. The prosecutors lack evidence and can only convict each prisoner for a lesser charge, so they offer the prisoners a deal. This deal comes as a dilemma to the prisoners as they need to choose between remaining silent or betraying their partner which would grant them a sentence lighter than that of the lesser charge, only if their partner remains silent. These two options can be represented as strategies in a game where remaining silent is referred to as cooperation and betraying the partner-in-crime is referred to as defection. This game can be applied to any situation in which there is a temptation for individuals to defect, however the net benefit of all parties is maximized if all individuals cooperate. It has been used to study the emergence of cooperation in a wide array of models in the fields of ecology, economics and psychology \cite{clark2001sequential,weitz2016oscillating,wong2005dynamic}.

In the PD, the payoff matrix for a given player is usually represented by:
\begin{equation}
    \renewcommand{\arraystretch}{1.45}
    \setlength\tabcolsep{3pt}
    \begin{tabular}{cc|c|c|}
      & \multicolumn{1}{c}{} & \multicolumn{2}{c}{Player $2$}\\
      & \multicolumn{1}{c}{} & \multicolumn{1}{c}{$D$}  & \multicolumn{1}{c}{$C$} \\\cline{3-4}
      \multirow{2}*{Player $1\:$}  & $D\:$ & $(a,a)$ & $(b,c)$ \\\cline{3-4}
      & $C\:$ & $(c,b)$ & $(d,d)$ \\\cline{3-4}
    \end{tabular}
\end{equation}
where every cell corresponds to each player choosing a strategy of either defect $(D)$ or cooperate $(C)$. The first and second elements of the tuple within each cell represent the payoff of players 1 and 2 respectively. To represent the dilemma, the payoffs are formulated with $b>d>a>c$ and to have the net payoff maximized for two cooperators, the system is further restricted such that $2d>b+c$. A common payoff matrix satisfying these conditions is:
\begin{equation}
    \setlength\tabcolsep{3pt}
    \renewcommand{\arraystretch}{1.45}
    \begin{tabular}{cc|c|c|}
      & \multicolumn{1}{c}{} & \multicolumn{2}{c}{Player $2$}\\
      & \multicolumn{1}{c}{} & \multicolumn{1}{c}{$D$}  & \multicolumn{1}{c}{$C$} \\\cline{3-4}
      \multirow{2}*{Player $1\:$}  & $D\:$ & $(1,1)$ & $(b,0)$ \\\cline{3-4}
      & $C\:$ & $(0,b)$ & $(3,3)$ \\\cline{3-4}
    \end{tabular}
    \label{eq:1}
\end{equation}
where $b>3$ is a parameter referred to as the temptation to defect. When simulated as a two-player game, the players' strategies will always converge to defection since it is the Nash equilibrium \cite{rubinstein1986finite}. However when iterated on a spatial structure with local interactions, more complex behaviour arises, including the persistence of cooperation due to the spatial clustering of alike strategies \cite{nowak1992evolutionary,nowak1994spatial}.

\section{Spatial Algebraic Model}
\label{sec:2}
For the model presented in this paper, the PD is iterated on rectangular lattice with  periodic boundary conditions where each cell represents a player with one of two strategies; defection represented by `0' and cooperation represented by `1'. A small $2\times 3$ lattice is used here due to current resource constraints of the computational algebraic analysis, but is illustrative of the general phenomena that arise.
\begin{SCfigure}[2.5][h]
    \hspace{0.5cm}
    \centering
    \resizebox{!}{0.1\textheight}{
    \begin{tikzpicture}
        \draw (0,0)--(0,2)--(3,2)--(3,0)--(0,0);
        \draw (1,0)--(1,2);
        \draw (2,0)--(2,2);
        \draw (0,1)--(3,1);
        \node at (0.5,1.5) {1};
        \node at (0.5,0.5) {2};
        \node at (1.5,1.5) {3};
        \node at (1.5,0.5) {4};
        \node at (2.5,1.5) {5};
        \node at (2.5,0.5) {6};
    \end{tikzpicture}
    }
    \caption{The spatial arrangement and enumeration of the cells on the $2\times 3$ spatial Prisoner's Dilemma lattice}
    \label{lattice}
\end{SCfigure}

\noindent The state space,
\begin{equation}
    X_{bin}=\{000000,000001,000010,000011,...,111110,111111\}
\end{equation}
is made of 64 6-bit binary strings, where the $i$\textsuperscript{th} bit from the left represents the strategy of cell $i$ (Figure \ref{lattice}). 
For a more notationally compact representation, this state set can also be written in decimal form with each state being the decimal integer equivalent of the binary string,
\begin{equation}
    X=\{0,1,2,3,...,62,63\}.
\end{equation}

During each synchronous playing of the game, $t$ (which we call a time step), each cell plays the PD with their von Neumann neighbours and gains a net payoff over all games using the payoff matrix (\ref{eq:1}). Note that since the system is a $2\times 3$ lattice, each cell has 3 von Neumann neighbours to avoid double-counting existing neighbours with the periodic boundaries. After playing against each other, each cell updates its strategy to match that of their neighbour with maximal payoff, only if the maximal payoff is greater than their own. If two neighbouring cells with different strategies have the same maximal payoff, then cooperation is chosen. 

To allow for more complexity, the model was formulated such that certain cells can have their strategy perturbed outside of $t$-dependent strategy evolution. We call these cells ``open''. If cell $i$ is open, there are two locally constant mappings associated with that cell; $d_i$ and $c_i$. These correspond to mapping the strategy of cell $i$ to defection or cooperation respectively, regardless of the change in payoff. These mappings on the set of open cells (denoted $O$) make a set of locally constant mappings, resetting cell $i$'s strategy to either $d$ or $c$ but leaving others' unchanged.
\begin{equation}
    T'_O=\{d_i,c_i\}_{i\in O}.
\end{equation}
The set of generators for the semigroup transformations is then given by
\begin{equation}
    T_O=T'_O\bigcup\{t\}.
\end{equation}
Words made from elements of $T_O$ define mappings on the set of states by applying each transformation in order from left to right. The set of transformations generated from $T_O$ comprise a semigroup denoted by
\begin{equation}
    S_O=\langle T_O \rangle,
\end{equation}
and $S_O$ acting on $X$ gives us the transformation semigroup $(X,S_O)$. As $T'_O$ is a set of locally constant mappings, it does not depend on the parameter $b$, however $t$ does and its $b$-dependence was explored using a \textit{python} script which also generated the semigroup mappings. In this analysis the strict inequality $b>3=d$ was relaxed so that $b\geq 3$. Note that for $b=3$, the system still favours defection since although mutual cooperation has become a weak Nash equilibrium, mutual defection is still the only strict Nash equilibrium, meaning no player can change their strategy without suffering a loss in payoff. The mappings generated by the \textit{python} script were then read into \textit{GAP} \cite{linton2007gap} and the transformation semigroup was analyzed using the \textit{SgpDec} package  \cite{egri2014sgpdec} to carry out a holonomy decomposition \cite{eilenberg76, nehaniv2015symmetry, egri2015computational}. This yields a KR decomposition of the spatial PD model's dynamics $(X,S_O)$ by identifying {\em natural subsystems}, i.e., nontrivial permutation groups whose state set is an image $X\cdot s$ of the state set $X$ under some semigroup element $s\in S_O$ and whose permutations are the restrictions of those members of $S_O$ which permute this set. Such an image set can be covered by the union of smaller image sets and singletons, which in turn must also be permuted by these transformations. The permutation group induced on the maximal covering sets of a natural subsystem by these sets is a {\em holonomy group}. See \cite{eilenberg76, nehaniv2015symmetry, egri2015computational} for details. In the next section, we will be referencing {\em subduction}, a generalized inclusion relation defined on the collection of images  together with $X$ and the singletons. For subsets $P,Q \subseteq X$, we say $P$ {\em subducts} $Q$ if $P \subseteq Q\cdot s$ for some $s \in S$ or $s$ the identity mapping. Mutual subduction implies isomorphism of holonomy groups, so equivalent locally reversible dynamics in the hierarchical decomposition can be compressed  \cite{egri2015computational}, giving insight into complexity of a dynamical system $(X,S)$. In the analysis and diagrams below, subduction corresponds to subset inclusion.                                
\section{Complexity Regimes}
\label{sec:3}

The investigation of the iterated PD's $b$-dependence revealed four different regimes characterized by unique sets of transformations by $t$ (Table \ref{regimetab}). The complexity of each regime was explored using the Krohn-Rhodes (KR) definition of semigroup complexity~\cite{krohn1968complexity}: KR complexity is formulated such that the complexity of a transformation semigroup $(X,S)$ is equal to the smallest number of non-trivial groups needed for a wreath product decomposition of $(X,S)$. Therefore an upper bound for the KR complexity is the number of levels with non-trivial groups in the holonomy decomposition. For the remainder of this paper, upper bounds will be used when referring to KR complexity.

\begin{SCtable}
\hspace{0.5cm}
\begin{tabular}{p{1.7cm}p{3.5cm}}
\hline\noalign{\smallskip}
\textbf{Regime} & \textbf{Parameter Range}  \\
\noalign{\smallskip}\Xhline{2\arrayrulewidth}\noalign{\smallskip}
A & $b>4.0$ \\ 
B & $b=4.0$ \\ 
C & $3.0<b<4.0$ \\
D & $b=3.0$ \\
\noalign{\smallskip}\hline\noalign{\smallskip}
\end{tabular}
\caption{Four unique regimes of the iterated Prisoner's Dilemma}
\label{regimetab} 
\end{SCtable}

\subsection{Regime A}

Beginning with regime A $(b>4.0)$, the system has a temptation to defect so large, that $t^2$ acting on any state containing at least one defector will bring that state to `000000' (state 0), which we will call pure defection.
 (As $t$ maps the pure cooperation state `111111' (state 63) to itself and no other states map to 63 by words generated by $t$, this state is left out of the subduction chains shown in Figures \ref{subA} and \ref{subB}.) The defection attractor dynamics can be visualized from subduction chain for $(X\setminus\{63\},\langle t \rangle)$ (Figure~\ref{subA}).  We can choose to only examine the mappings induced by words generated by $t$ when comparing regimes since the semigroup generated by $T'_O$ is unchanged by the parameter $b$. As $t^2$ only maps to pure defection and the rest of the mappings in $S_O$ are locally constant maps, there are no pools of reversibility and few levels in the holonomy decomposition, yielding a relatively trivial system. 

\begin{figure}[h]
    \centering
    \begin{minipage}{0.4\textwidth}
    \resizebox{!}{0.14\textheight}{
    \begin{tikzpicture}[node distance=1.3cm]
        \node (top) [subpos] {$X\setminus\{63\}$};
        \node (mid) [subpos, below of=top] {$\{0,1,2,4,5,8,10,16,17,20,32,34,40\}$};
        \node (bot) [subpos, below of=mid] {$\{0\}$};
        \draw [arrow] (top)--(mid);
        \draw [arrow] (mid)--(bot);
    \end{tikzpicture}
    }
    \caption{Subduction chain for $(X,\langle t \rangle)$ with $b>4$}
    \label{subA}
    \end{minipage}\hfill
    \begin{minipage}{0.55\textwidth}
    \resizebox{!}{0.14\textheight}{
    \begin{tikzpicture}[node distance=1.3cm]
        \node (top) [subpos] {$X\setminus\{63\}$};
        \node (mid) [subpos, below of=top] {$\{0,5,10,17,20,21,23,29,34,40,42,43,46,53,58\}$};
        \node (bot) [subpos, below of=mid] {$\{0,21,23,29,42,43,46,53,58\}$};
        \draw [arrow] (top)--(mid);
        \draw [arrow] (mid)--(bot);
    \end{tikzpicture}
    }
    \caption{Subduction chain for $(X,\langle t \rangle)$ with $b=4$}
    \label{subB}
    \end{minipage}
\end{figure}

\subsection{Regime B}
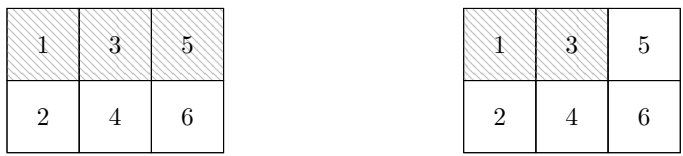
\begin{figure}
    \centering
    \begin{tabular}{c @{\quad} c }
    \resizebox{!}{0.1\textheight}{
    \begin{tikzpicture}
        \draw (0,0)--(0,2)--(3,2)--(3,0)--(0,0);
        \draw (1,0)--(1,2);
        \draw (2,0)--(2,2);
        \draw (0,1)--(3,1);
        \draw [pattern=north west lines, pattern color=gray!50] (0,2) rectangle (1,1);
        \draw (0,1) rectangle (1,0);
        \draw [pattern=north west lines, pattern color=gray!50] (1,2) rectangle (2,1);
        \draw (1,1) rectangle (2,0);
        \draw [pattern=north west lines, pattern color=gray!50] (2,2) rectangle (3,1);
        \draw (2,1) rectangle (3,0);
        \node at (0.5,1.5) {1};
        \node at (0.5,0.5) {2};
        \node at (1.5,1.5) {3};
        \node at (1.5,0.5) {4};
        \node at (2.5,1.5) {5};
        \node at (2.5,0.5) {6};
    \end{tikzpicture}
    } &
    \resizebox{!}{0.1\textheight}{
    \begin{tikzpicture}
        \draw (0,0)--(0,2)--(3,2)--(3,0)--(0,0);
        \draw (1,0)--(1,2);
        \draw (2,0)--(2,2);
        \draw (0,1)--(3,1);
        \draw [pattern=north west lines, pattern color=gray!50] (0,2) rectangle (1,1);
        \draw (0,1) rectangle (1,0);
        \draw [pattern=north west lines, pattern color=gray!50] (1,2) rectangle (2,1);
        \draw (1,1) rectangle (2,0);
        \draw (2,2) rectangle (3,1);
        \draw (2,1) rectangle (3,0);
        \node at (0.5,1.5) {1};
        \node at (0.5,0.5) {2};
        \node at (1.5,1.5) {3};
        \node at (1.5,0.5) {4};
        \node at (2.5,1.5) {5};
        \node at (2.5,0.5) {6};
    \end{tikzpicture}
    } \\[2mm]
    \footnotesize (a) ``3-in-a-row'' strategy configuration &
    \footnotesize(b) ``L-shape'' strategy configuration represented \\
    \footnotesize represented by $\{21,42\}=[21]_{\cong} \subset X$ & 
    \footnotesize by $\{23, 29, 43, 46, 53, 58\}=[23]_{\cong}\subset X$
  \end{tabular}
  \caption{Mixed strategy equilibria configuration for regimes B and C. Hatched pattern and no fill represent defector and cooperator strategies, respectively.}
    \label{mixeq}
\end{figure}

Regime B $(b=4.0)$ can be seen as a critical point where the system changes from regime A to C. The main difference between regimes A and B is that mixed strategy equilibria under transformation $t$ appear in regime B. These equilibria fall under two spatial configurations up to isomorphism: ``3-in-a-row’’ and ``L-shape’’, shown in Figure \ref{mixeq}. Similar to regime A, this regime does not have non-trivial groups in the holonomy decomposition giving both regimes a KR complexity of 0.

\subsection{Regime C}
\begin{SCfigure}[2][h]
    \hspace{0.5cm}
    \centering
    \resizebox{!}{0.14\textheight}{
    \begin{tikzpicture}[node distance=1.2cm]
        \node (top) [subpos] {$X$};
        \node (mid) [subpos, below of=top] {$\{0,5,10,17,20,23,29,34,40,43,46,53,58,63\}$};
        \node (bot) [subpos, below of=mid] {$\{0,23,29,43,46,53,58,63\}$};
        \draw [arrow] (top)--(mid);
        \draw [arrow] (mid)--(bot);
    \end{tikzpicture}
    }
    \caption{Subduction chain for $(X,\langle t\rangle)$ with $3<b<4$}
    \label{subC}
\end{SCfigure}
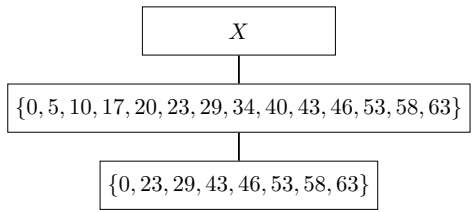

In regime C $(3.0<b<4.0)$, the temptation to defect is at an intermediate level which now allows certain states to map to ones of higher cooperation with $t$. From the subduction chain (Figure \ref{subC}) one can see that the decreased temptation to defect removes one class of mixed strategy equilibria states, $[21]_{\cong}$. In this state, the defector's net payoff is $b+2$ as it receives $b$ for playing against one adjacent cooperator as well as 2 for playing against two adjacent defectors. Since $b<4$ in this regime, $b+2<6$, the total payoff for cooperators and this state will now map to pure cooperation when acted on by $t$.

This regime is drastically different from regimes A and B as there are now cyclic groups in the holonomy decomposition. For these intermediate temptations to defect, the system has pools of reversibility where dynamical cycles may recur, unlike the previous regimes where any non-trivial mappings induced by words in $S_O$ would bring the system to a state in which the previous state is inaccessible by the same transformation. This reversibility is entirely $t$-dependent since any words made of locally constant maps which act non-trivially on a state are by definition irreversible. Only when the temptation to defect is low enough such that an action by $t$ can bring the system to a new state of equivalent or higher cooperation will mappings induced by words from $S_O$ form non-trivial groups.

Additionally, the amount and distribution of open cells now play a significant role in the system complexity. In general, KR complexity increases with the number of open cells; yet for a given number of open cells, their distribution plays a significant role (Table \ref{opencells}). Note the open cell configuration corresponding to $T'_{123}=\{d_1,d_2,d_3,c_1,c_2,c_3\}$ has its upper bound of KR complexity reduced from 6 to 4. This is because the configuration $T'_{1234}$ has KR complexity 4 and $T'_{123}$ is a subsemigroup of $T'_{1234}$. It follows naturally that $(X,S_{1234})$ emulates $(X,S_{123})$ and from the KR complexity axioms \cite{krohn1968complexity}, for  transformation semigroups $(Y,T)$ and $(X,S)$, if $(X,S)$ can emulate $(Y,T)$ then the complexity of $(Y,T)$ must be less than or equal to that of $(X,S)$. 

\begin{table}[!hbt]
\caption{Unique open cell configurations for 2,3 \& 4 open cells, represented in grey}
\label{opencells}
\resizebox{\columnwidth}{!}{
\begin{tabular}{p{0.173\linewidth}p{0.084\linewidth}p{0.084\linewidth}p{0.084\linewidth}p{0.084\linewidth}p{0.084\linewidth}p{0.084\linewidth}p{0.084\linewidth}p{0.084\linewidth}p{0.084\linewidth}}
\hline\noalign{\smallskip}
    \textbf{\begin{tabular}[c]{@{}l@{}}Open Cell\\ Orientation\end{tabular}}
    &
    \resizebox{!}{0.03\textheight}{
    \begin{tikzpicture}
        \draw (0,0)--(0,2)--(3,2)--(3,0)--(0,0);
        \draw (1,0)--(1,2);
        \draw (2,0)--(2,2);
        \draw (0,1)--(3,1);
        \draw [fill=gray!50] (0,2) rectangle (1,1);
        \draw [fill=gray!50] (1,2) rectangle (2,1);
    \end{tikzpicture}
    }
    &
    \resizebox{!}{0.03\textheight}{
     \begin{tikzpicture}
        \draw (0,0)--(0,2)--(3,2)--(3,0)--(0,0);
        \draw (1,0)--(1,2);
        \draw (2,0)--(2,2);
        \draw (0,1)--(3,1);
        \draw [fill=gray!50] (0,2) rectangle (1,1);
        \draw [fill=gray!50] (1,1) rectangle (2,0);
    \end{tikzpicture}
    }
    &
    \resizebox{!}{0.03\textheight}{
     \begin{tikzpicture}
        \draw (0,0)--(0,2)--(3,2)--(3,0)--(0,0);
        \draw (1,0)--(1,2);
        \draw (2,0)--(2,2);
        \draw (0,1)--(3,1);
        \draw [fill=gray!50] (0,2) rectangle (1,1);
        \draw [fill=gray!50] (0,1) rectangle (1,0);
    \end{tikzpicture}
    }
    &
    \resizebox{!}{0.03\textheight}{
     \begin{tikzpicture}
        \draw (0,0)--(0,2)--(3,2)--(3,0)--(0,0);
        \draw (1,0)--(1,2);
        \draw (2,0)--(2,2);
        \draw (0,1)--(3,1);
        \draw [fill=gray!50] (0,2) rectangle (1,1);
        \draw [fill=gray!50] (1,1) rectangle (2,0);
        \draw [fill=gray!50] (2,2) rectangle (3,1);
    \end{tikzpicture}
    }
    & 
    \resizebox{!}{0.03\textheight}{
    \begin{tikzpicture}
        \draw (0,0)--(0,2)--(3,2)--(3,0)--(0,0);
        \draw (1,0)--(1,2);
        \draw (2,0)--(2,2);
        \draw (0,1)--(3,1);
        \draw [fill=gray!50] (0,2) rectangle (1,1);
        \draw [fill=gray!50] (1,2) rectangle (2,1);
        \draw [fill=gray!50] (2,2) rectangle (3,1);
    \end{tikzpicture}
    }
    &
    \resizebox{!}{0.03\textheight}{
     \begin{tikzpicture}
        \draw (0,0)--(0,2)--(3,2)--(3,0)--(0,0);
        \draw (1,0)--(1,2);
        \draw (2,0)--(2,2);
        \draw (0,1)--(3,1);
        \draw [fill=gray!50] (0,2) rectangle (1,1);
        \draw [fill=gray!50] (0,1) rectangle (1,0);
        \draw [fill=gray!50] (1,2) rectangle (2,1);
    \end{tikzpicture}
    }
    &
    \resizebox{!}{0.03\textheight}{
     \begin{tikzpicture}
        \draw (0,0)--(0,2)--(3,2)--(3,0)--(0,0);
        \draw (1,0)--(1,2);
        \draw (2,0)--(2,2);
        \draw (0,1)--(3,1);
        \draw [fill=gray!50] (0,2) rectangle (1,1);
        \draw [fill=gray!50] (0,1) rectangle (1,0);
        \draw [fill=gray!50] (1,2) rectangle (2,1);
        \draw [fill=gray!50] (1,1) rectangle (2,0);
    \end{tikzpicture}
    }
    &
    \resizebox{!}{0.03\textheight}{
     \begin{tikzpicture}
        \draw (0,0)--(0,2)--(3,2)--(3,0)--(0,0);
        \draw (1,0)--(1,2);
        \draw (2,0)--(2,2);
        \draw (0,1)--(3,1);
        \draw [fill=gray!50] (0,2) rectangle (1,1);
        \draw [fill=gray!50] (1,2) rectangle (2,1);
        \draw [fill=gray!50] (1,1) rectangle (2,0);
        \draw [fill=gray!50] (2,1) rectangle (3,0);
    \end{tikzpicture}
    }
    &
    \resizebox{!}{0.03\textheight}{
     \begin{tikzpicture}
        \draw (0,0)--(0,2)--(3,2)--(3,0)--(0,0);
        \draw (1,0)--(1,2);
        \draw (2,0)--(2,2);
        \draw (0,1)--(3,1);
        \draw [fill=gray!50] (0,2) rectangle (1,1);
        \draw [fill=gray!50] (0,1) rectangle (1,0);
        \draw [fill=gray!50] (1,2) rectangle (2,1);
        \draw [fill=gray!50] (2,2) rectangle (3,1);
    \end{tikzpicture}
    }\\
    \noalign{\smallskip}\Xhline{2\arrayrulewidth}\noalign{\smallskip}
    \textbf{\begin{tabular}[c]{@{}l@{}}KR\\Complexity\end{tabular}} & 0 & 0 & 2 & 0 & 2 & $\cancel{6}$ 4 & 4 & 4 & 7\\
    \textbf{\begin{tabular}[c]{@{}l@{}}Groups in\\Holonomy\\Decomposition\end{tabular}} & --- & --- & \begin{tabular}[c]{@{}l@{}} $(3,C_2)$ \\ $ (2,C_2)$\end{tabular} & --- & $(2,C_2)$ & \begin{tabular}[c]{@{}l@{}} $(4,C_2)$ \\ $(3,C_2)$ \\ $ (2,C_2)$ \end{tabular} & \begin{tabular}[c]{@{}l@{}} $(4,C_2)$ \\ $(3,C_2)$ \\ $ (2,C_2)$ \end{tabular} & \begin{tabular}[c]{@{}l@{}} $(4,C_2)$ \\ $(3,C_2)$ \\ $ (2,C_2)$ \end{tabular} & \begin{tabular}[c]{@{}l@{}} $(3,S_3)$ \\ $(4,C_2)$ \\ $(3,C_2)$ \\ $ (2,C_2)$ \end{tabular}\\
    \noalign{\smallskip}\hline\noalign{\smallskip}
\end{tabular}
}
\end{table}

\begin{figure}[h]
    \vspace{-0.5cm}
    \centering
    \begin{tabular}{c @{} c }
    \includegraphics[width=40mm,trim={13mm 5mm 10mm 0},clip]{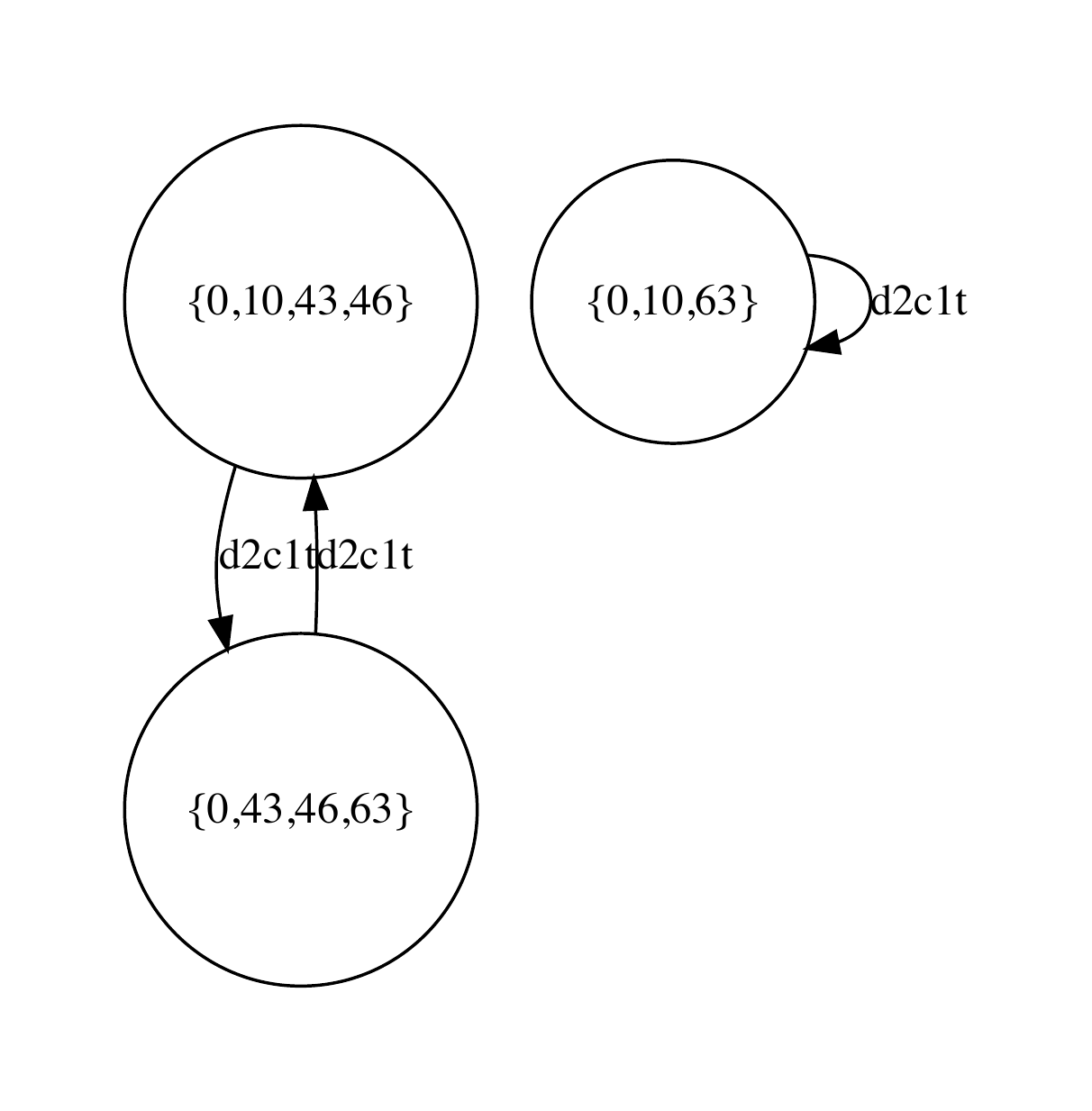}
    &
    \includegraphics[width=75mm,trim={10mm 5mm 0 0},clip]{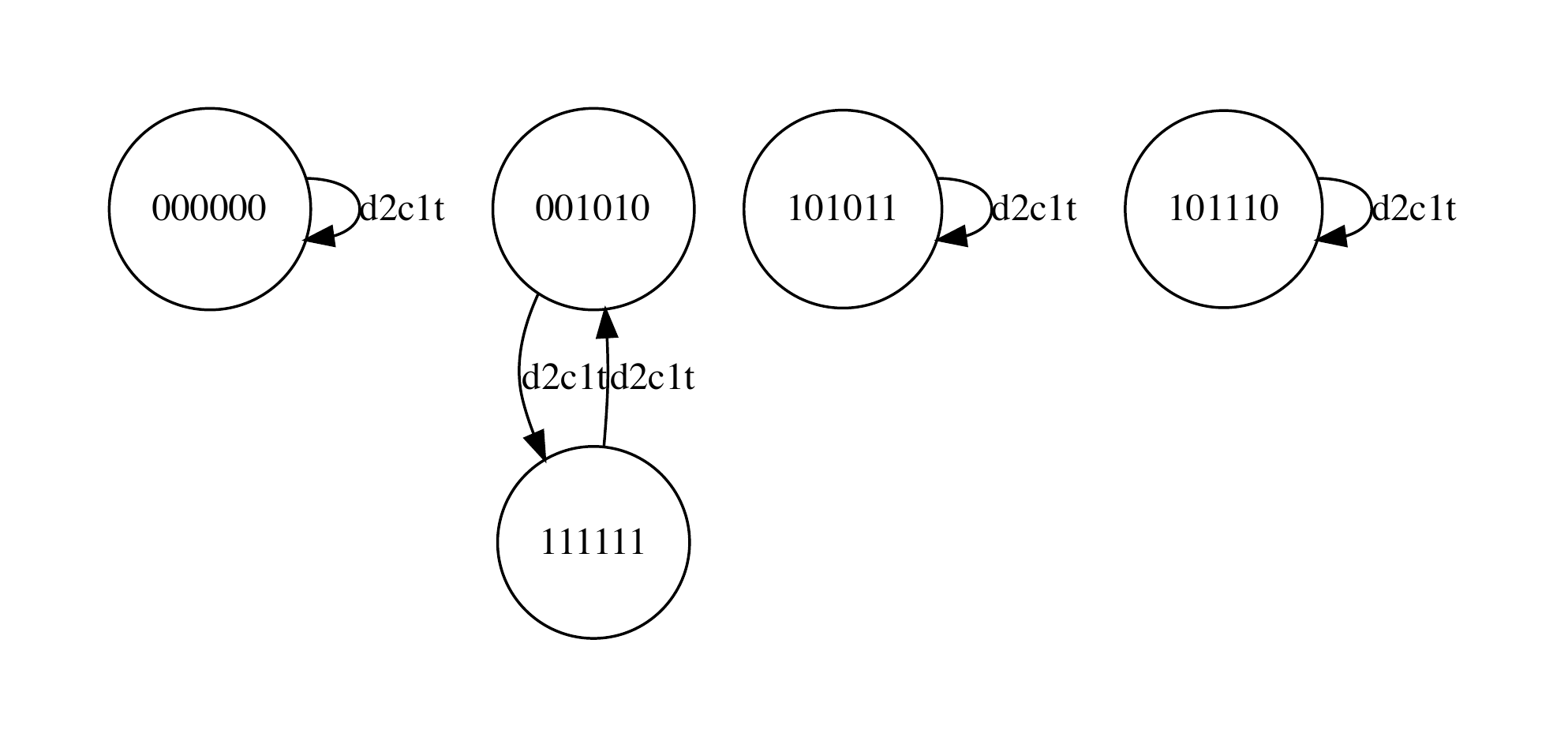}\\
     \footnotesize (a) Orbits of holonomy group $(3,C_2)$  & \footnotesize (b) Natural subsystem for transformation $d_2c_1t$\\ 
     \footnotesize for transformation $d_2c_1t$ &
  \end{tabular}
    \caption{Orbits and natural subsystem for a $(3,C_2)$ found in regime $C$ with $O=\{1,2\}$. Note states numbered 0, 10, 43, 45, and 63 in (a) correspond to strategy distributions  `000000', `001010', `101011', `101110', and `111111' seen in (b), respectively.}
    \label{d2c1t}
\end{figure}

\noindent Below are the orbits for a representation of the holonomy group $(3,C_2)$ found in the holonomy decomposition (Figure \ref{d2c1t}a). This is one of two $C_2$ permutator groups in regime C with $O=\{1,2\}$. The generator of this permutation group is $d_2c_1t$ which represents mapping cell 2 to defection, cell 1 to cooperation and then letting one time step, $t$ pass.
Since the holonomy group $(3,C_2)$ shows the group action on a set of 3 subsets permuted by permutations of 5 underlying states, the exact mechanism leading to this reversibility is not immediately clear. We can gain a better understanding of the dynamics of this cyclic group by examining the orbits of the transformation $d_2c_1t$ on specific states in these sets as shown in the natural subsystem (Figure \ref{d2c1t}b). 

Most of these orbits act in an expected manner since for the two right-most orbits, the locally constant mappings do not change the state and as the states are mixed strategy equilibria for regime C, action by $t$ does not change the state. In the left-most orbit, the behaviour is also expected since $d_2c_1$ effectively turns a pure defection state into one with a single cooperator which will receive the lowest payoff of its defecting neighbours, thus reverting back to defection with $t$. 

In the orbit second from the left in Fig.~\ref{d2c1t}b, the behaviour is much more interesting as the same transformation that removes four cooperators from the system, also leads it into a state of pure cooperation. For state `001010', $d_2c_1$ acts as simply $c_1$ since cell 2 is already a defector. This maps the system to state `101010'. As shown above, this state now maps to one of pure cooperation with $t$. At the state of pure cooperation, $d_2c_1$ now acts as $d_2$ mapping the state to `101111'.

\begin{figure}[h]
\centering
\begin{tabular}{>{\centering\arraybackslash}m{26mm}>{\centering\arraybackslash}m{9mm}>{\centering\arraybackslash}m{26mm}>{\centering\arraybackslash}m{9mm}>{\centering\arraybackslash}m{26mm}>{\centering\arraybackslash}m{9mm}}
    \centering
    \resizebox{!}{0.09\textheight}{
    \begin{tikzpicture}
        \draw (0,0)--(0,2)--(3,2)--(3,0)--(0,0);
        \draw (1,0)--(1,2);
        \draw (2,0)--(2,2);
        \draw (0,1)--(3,1);
        \draw [pattern=north west lines, pattern color=gray!50] (0,2) rectangle (1,1);
        \draw [pattern=north west lines, pattern color=gray!50] (0,1) rectangle (1,0);
        \draw (1,2) rectangle (2,1);
        \draw [pattern=north west lines, pattern color=gray!50] (1,1) rectangle (2,0);
        \draw (2,2) rectangle (3,1);
        \draw [pattern=north west lines, pattern color=gray!50] (2,1) rectangle (3,0);
        \node at (0.5,1.5) {1};
        \node at (0.5,0.5) {2};
        \node at (1.5,1.5) {3};
        \node at (1.5,0.5) {4};
        \node at (2.5,1.5) {5};
        \node at (2.5,0.5) {6};
    \end{tikzpicture}
    }
&
$\xrightarrow{\: d_2c_1 \:}$
&
    \resizebox{!}{0.09\textheight}{
    \begin{tikzpicture}
        \draw (0,0)--(0,2)--(3,2)--(3,0)--(0,0);
        \draw (1,0)--(1,2);
        \draw (2,0)--(2,2);
        \draw (0,1)--(3,1);
        \draw (0,2) rectangle (1,1);
        \draw [pattern=north west lines, pattern color=gray!50] (0,1) rectangle (1,0);
        \draw (1,2) rectangle (2,1);
        \draw [pattern=north west lines, pattern color=gray!50] (1,1) rectangle (2,0);
        \draw (2,2) rectangle (3,1);
        \draw [pattern=north west lines, pattern color=gray!50] (2,1) rectangle (3,0);
        \node at (0.5,1.5) {1};
        \node at (0.5,0.5) {2};
        \node at (1.5,1.5) {3};
        \node at (1.5,0.5) {4};
        \node at (2.5,1.5) {5};
        \node at (2.5,0.5) {6};
    \end{tikzpicture}
    }
&
$\xrightarrow{\quad t \quad}$
&
    \resizebox{!}{0.09\textheight}{
    \begin{tikzpicture}
        \draw (0,0)--(0,2)--(3,2)--(3,0)--(0,0);
        \draw (1,0)--(1,2);
        \draw (2,0)--(2,2);
        \draw (0,1)--(3,1);
        \draw (0,2) rectangle (1,1);
        \draw (0,1) rectangle (1,0);
        \draw (1,2) rectangle (2,1);
        \draw (1,1) rectangle (2,0);
        \draw (2,2) rectangle (3,1);
        \draw (2,1) rectangle (3,0);
        \node at (0.5,1.5) {1};
        \node at (0.5,0.5) {2};
        \node at (1.5,1.5) {3};
        \node at (1.5,0.5) {4};
        \node at (2.5,1.5) {5};
        \node at (2.5,0.5) {6};
    \end{tikzpicture}
    }\\
    001010 & & 101010 & & 111111\\
    \noalign{\vskip 3mm}
    \centering
    \resizebox{!}{0.09\textheight}{
    \begin{tikzpicture}
        \draw (0,0)--(0,2)--(3,2)--(3,0)--(0,0);
        \draw (1,0)--(1,2);
        \draw (2,0)--(2,2);
        \draw (0,1)--(3,1);
        \draw (0,2) rectangle (1,1);
        \draw (0,1) rectangle (1,0);
        \draw (1,2) rectangle (2,1);
        \draw (1,1) rectangle (2,0);
        \draw (2,2) rectangle (3,1);
        \draw (2,1) rectangle (3,0);
        \node at (0.5,1.5) {1};
        \node at (0.5,0.5) {2};
        \node at (1.5,1.5) {3};
        \node at (1.5,0.5) {4};
        \node at (2.5,1.5) {5};
        \node at (2.5,0.5) {6};
    \end{tikzpicture}
    }
&
$\xrightarrow{\: d_2c_1 \:}$
&
    \resizebox{!}{0.09\textheight}{
    \begin{tikzpicture}
        \draw (0,0)--(0,2)--(3,2)--(3,0)--(0,0);
        \draw (1,0)--(1,2);
        \draw (2,0)--(2,2);
        \draw (0,1)--(3,1);
        \draw (0,2) rectangle (1,1);
        \draw [pattern=north west lines, pattern color=gray!50] (0,1) rectangle (1,0);
        \draw (1,2) rectangle (2,1);
        \draw (1,1) rectangle (2,0);
        \draw (2,2) rectangle (3,1);
        \draw (2,1) rectangle (3,0);
        \node at (0.5,1.5) {1};
        \node at (0.5,0.5) {2};
        \node at (1.5,1.5) {3};
        \node at (1.5,0.5) {4};
        \node at (2.5,1.5) {5};
        \node at (2.5,0.5) {6};
    \end{tikzpicture}
    }
&
$\xrightarrow{\quad t \quad}$
&
    \resizebox{!}{0.09\textheight}{
    \begin{tikzpicture}
        \draw (0,0)--(0,2)--(3,2)--(3,0)--(0,0);
        \draw (1,0)--(1,2);
        \draw (2,0)--(2,2);
        \draw (0,1)--(3,1);
        \draw [pattern=north west lines, pattern color=gray!50] (0,2) rectangle (1,1);
        \draw [pattern=north west lines, pattern color=gray!50] (0,1) rectangle (1,0);
        \draw (1,2) rectangle (2,1);
        \draw [pattern=north west lines, pattern color=gray!50] (1,1) rectangle (2,0);
        \draw (2,2) rectangle (3,1);
        \draw [pattern=north west lines, pattern color=gray!50] (2,1) rectangle (3,0);
        \node at (0.5,1.5) {1};
        \node at (0.5,0.5) {2};
        \node at (1.5,1.5) {3};
        \node at (1.5,0.5) {4};
        \node at (2.5,1.5) {5};
        \node at (2.5,0.5) {6};
    \end{tikzpicture}
    }\\
    111111 & & 101111 & & 001010 \\
\end{tabular}
\caption{Dynamics of $C_2$ generated by $d_2c_1t$ \label{fig7}}
\end{figure}
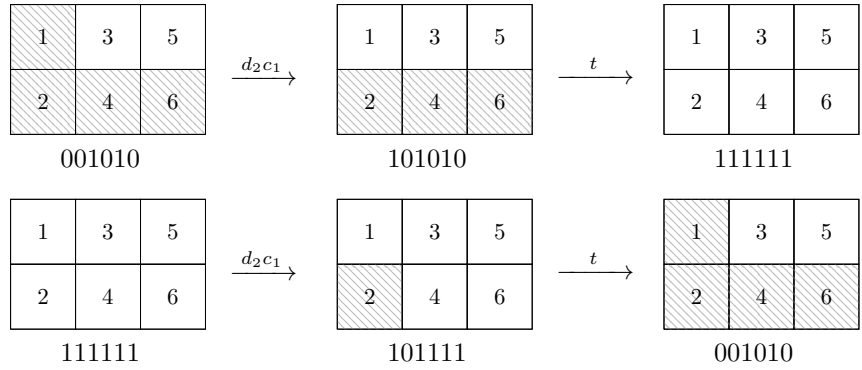

\begin{figure}[h]
    \vspace{-0.8cm}
    \centering
    \includegraphics[width=100mm,trim={-5mm 10mm 13mm 0},clip]{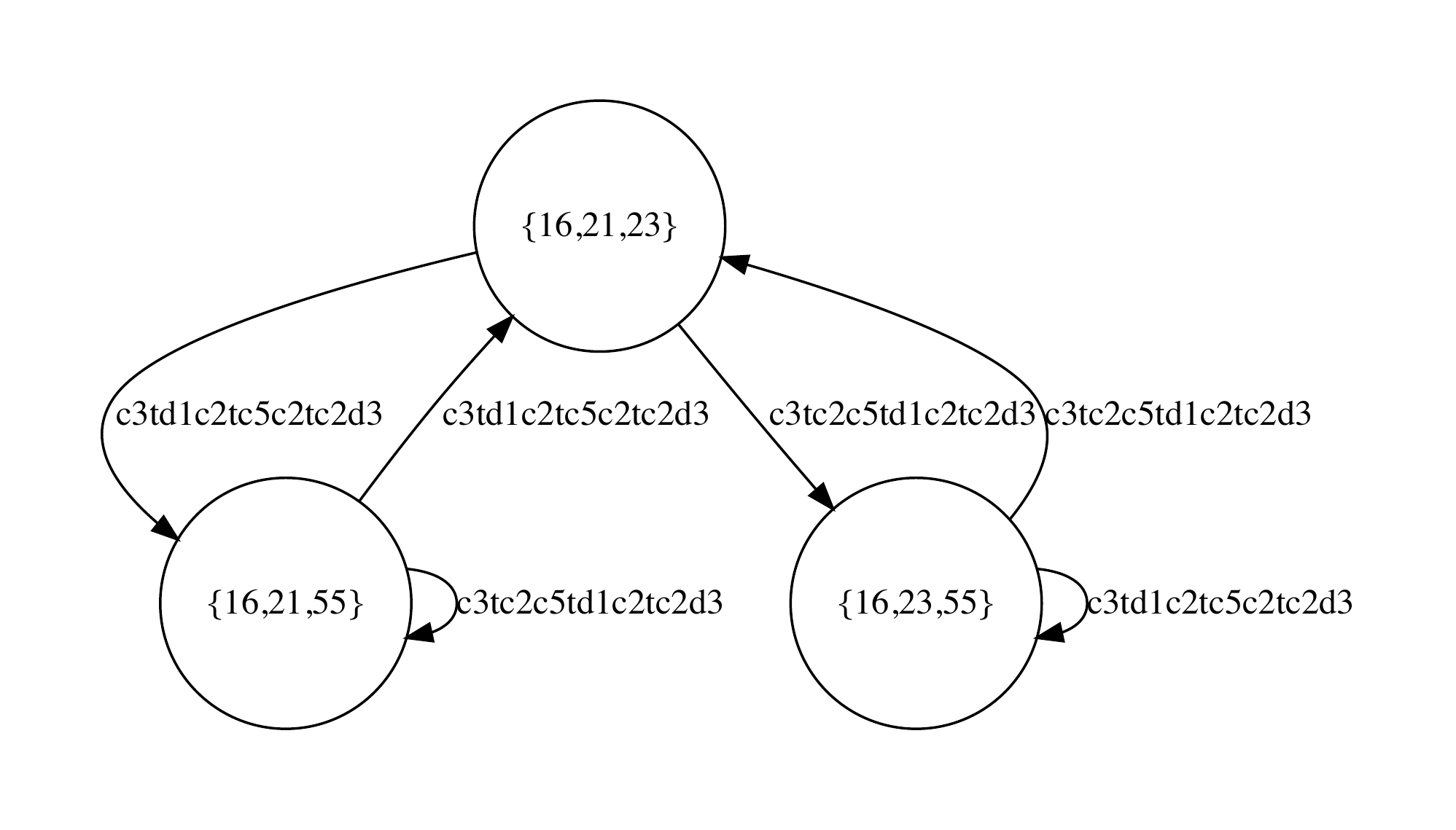}
    \vspace{-0.3cm}
    \caption{Orbits of Holonomy Group $(3,S_3)$. The underlying states (strategy distributions) of the corresponding natural subsystem are 16=`01000', 21=`010101', 23=`010111', and 55=`101101'.}
    \label{S3}
\end{figure}

\noindent In this state, the single defector benefits from being surrounded by cooperators receiving the highest net payoff as $b>3$ and with $t$, all of its neighbours switch to the defector strategy resulting in the state `001010' (See Figure~\ref{fig7}). 
These pools of reversibility can offer some insight regarding the spatial configuration of strategies which lead to the persistence of cooperation.
Additionally, in regime C, we see the only symmetric group $(3,S_3)$ in the holonomy decomposition (Figure \ref{S3}). Here the two group generators are given by $c_3td_1c_2tc_5c_2tc_2d_3$ and $c_3tc_2c_5td_1c_2tc_2d_3$  and although these words are long and hard to interpret, the possibility of appearance of such non-abelian group dynamics is not an obvious result of simple iterated PD.

\subsection{Regime D}

In regime D $(b=3.0)$, which can be interpreted as a weak PD, the temptation to defect is very low and consequently the incentive to cooperate is highest. Due to this push towards cooperation, there are less pools of reversibility than in regime C and the highest upper bound for KR complexity is 2. From the subduction chain (Figure~\ref{subD}), we see that all equilibria are mapped to by a single time step $t$ and a new class of equilibria emerge. (Also present are the ``L-shape" equilibria with two defectors in a row we encountered above.)  This new class represents a single defector, which in all previous regimes had been beneficial to the lone defector. In this regime, the temptation to defect is low enough that the system has become immune to invasion by a single defector.

\begin{figure}[h]
    \centering
    \begin{minipage}{0.45\textwidth}
    \resizebox{!}{0.09\textheight}{
    \begin{tikzpicture}[node distance=0.1cm]
        \node (top) [subpos] {$X$};
        \node (bot) [subpos, below of=mid] {$\{0,23,29,31,43,46,47,53,55,58,59,61,62,63\}$};
        \draw [arrow] (top)--(bot);
    \end{tikzpicture}
    }
    \caption{Subduction chain for $(X,\langle t \rangle)$ with $b=3$}
    \label{subD}
    \end{minipage}\hfill
    \begin{minipage}{0.45\textwidth}
    \centering
    \resizebox{!}{0.1\textheight}{
    \begin{tikzpicture}
        \draw (0,0)--(0,2)--(3,2)--(3,0)--(0,0);
        \draw (1,0)--(1,2);
        \draw (2,0)--(2,2);
        \draw (0,1)--(3,1);
        \draw (0,2) rectangle (1,1);
        \draw (0,1) rectangle (1,0);
        \draw (1,2) rectangle (2,1);
        \draw (1,1) rectangle (2,0);
        \draw (2,2) rectangle (3,1);
        \draw [pattern=north west lines, pattern color=gray!50] (2,1) rectangle (3,0);
        \node at (0.5,1.5) {1};
        \node at (0.5,0.5) {2};
        \node at (1.5,1.5) {3};
        \node at (1.5,0.5) {4};
        \node at (2.5,1.5) {5};
        \node at (2.5,0.5) {6};
    \end{tikzpicture}
    }
    \caption{New class of equilibria $\{31,47,55,59,61,62\}=[31]_{\cong}\subset X$}
    \label{eqD}
    \end{minipage}
\end{figure}

\vspace{-0.4cm}

\section{Conclusion}

Representing the iterated PD as a transformation semigroup allows the holonomy decomposition to reveal qualitative differences between distinct payoff-dependent regimes. When the temptation to defect is below a threshold, the KR complexity becomes non-zero and pools of reversibility form. The number of open cells also positively influences the KR complexity, however their spatial distribution plays an equally important role.
With greater computational power, it would be interesting to further explore this system with a larger number of players as well as different topologies to see how the results presented in this paper compare to larger and more complex spatial configurations. With this information, one could explore how the KR complexity varies with both spatial size and configuration, as well as with the temptation to defect. Additionally, it could lead to insights resulting in algebraic theorems for more general iterated PD systems.
 
\pagebreak

\bibliographystyle{plain}
{\footnotesize\bibliography{mybib.bib}}

\begin{thebibliography}{10}

\bibitem{clark2001sequential}
Kenneth Clark and Martin Sefton.
\newblock The sequential prisoner's dilemma: evidence on reciprocation.
\newblock {\em The Economic Journal}, 111(468):51--68, 2001.

\bibitem{egri2014sgpdec}
Attila Egri-Nagy, James~D Mitchell, and Chrystopher~L Nehaniv.
\newblock Sgp{D}ec: Cascade (de)compositions of finite transformation
  semigroups and permutation groups.
\newblock In {\em International Congress on Mathematical Software}, pages
  75--82. Springer Lecture Notes in Computer Science, vol.\ 8592, 2014.

\bibitem{egri2015computational}
Attila Egri-Nagy and Chrystopher~L Nehaniv.
\newblock Computational holonomy decomposition of transformation semigroups.
\newblock {\em arXiv preprint arXiv:1508.06345}, 2015.

\bibitem{eilenberg76}
Samuel Eilenberg.
\newblock {\em Automata, Languages and Machines}, volume~B.
\newblock Academic Press, 1976.

\bibitem{krohn1965algebraic}
Kenneth Krohn and John Rhodes.
\newblock Algebraic theory of machines. {I}. prime decomposition theorem for
  finite semigroups and machines.
\newblock {\em Trans. American Mathematical Society}, 116:450--464, 1965.

\bibitem{krohn1968complexity}
Kenneth Krohn and John Rhodes.
\newblock Complexity of finite semigroups.
\newblock {\em Annals Math.}, pages 128--160, 1968.

\bibitem{linton2007gap}
Steve Linton.
\newblock {G}{A}{P}: groups, algorithms, programming.
\newblock {\em ACM Communications in Computer Algebra}, 41(3):108--109, 2007.

\bibitem{nehaniv2015symmetry}
Chrystopher~L Nehaniv, John Rhodes, Attila Egri-Nagy, Paolo Dini,
  Eric~Rothstein Morris, G{\'a}bor Horv{\'a}th, Fariba Karimi, Daniel
  Schreckling, and Maria~J Schilstra.
\newblock Symmetry structure in discrete models of biochemical systems: natural
  subsystems and the weak control hierarchy in a new model of computation
  driven by interactions.
\newblock {\em Philosophical Transactions of the Royal Society A: Mathematical,
  Physical and Engineering Sciences}, 373(2046):20140223, 2015.

\bibitem{nowak1994spatial}
Martin~A Nowak, Sebastian Bonhoeffer, and Robert~M May.
\newblock Spatial games and the maintenance of cooperation.
\newblock {\em Proceedings of the National Academy of Sciences},
  91(11):4877--4881, 1994.

\bibitem{nowak1992evolutionary}
Martin~A Nowak and Robert~M May.
\newblock Evolutionary games and spatial chaos.
\newblock {\em Nature}, 359(6398):826, 1992.

\bibitem{rubinstein1986finite}
Ariel Rubinstein.
\newblock Finite automata play the repeated prisoner's dilemma.
\newblock {\em Journal of Economic Theory}, 39(1):83--96, 1986.

\bibitem{weitz2016oscillating}
Joshua~S Weitz, Ceyhun Eksin, Keith Paarporn, Sam~P Brown, and William~C
  Ratcliff.
\newblock An oscillating tragedy of the commons in replicator dynamics with
  game-environment feedback.
\newblock {\em Proceedings of the National Academy of Sciences},
  113(47):E7518--E7525, 2016.

\bibitem{wong2005dynamic}
Rosanna Yin-mei Wong and Ying-yi Hong.
\newblock Dynamic influences of culture on cooperation in the prisoner's
  dilemma.
\newblock {\em Psychological science}, 16(6):429--434, 2005.

\end{thebibliography}

\fontsize{8pt}{12pt}\selectfont
\noindent All code used to generate semigroup mappings, analysis and figures is available at \href{https://gist.github.com/cello-kabob/eb0a56bbcd598bf613d91b8e773ff9cf}{\underline{https://git.io/JJJcN}}

\end{document}